\documentclass[preprint,12pt]{elsarticle}
\usepackage{amsmath}
\usepackage{amsfonts}
\usepackage{bm}
\usepackage{lipsum}
\usepackage{enumerate}
\usepackage{graphicx}
\usepackage{amsfonts}

\usepackage{mathrsfs}
\usepackage{subfigure}
\usepackage{epstopdf}
\usepackage{url}
\usepackage{verbatim}
\usepackage{amssymb}
\usepackage{hyperref}
\usepackage{color}
\usepackage{listings}

\def\@begintheorem#1#2{\par\bgroup{\scshape #1\ #2. }\it\ignorespaces}
\def\@opargbegintheorem#1#2#3{\par\bgroup%
   {\scshape #1\ #2\ ({\upshape #3}). }\it\ignorespaces}
\def\@endtheorem{\egroup}

\if@onethmnum
  \newtheorem{theorem}{Theorem}
  \newtheorem{lemma}[theorem]{Lemma}
  \newtheorem{corollary}[theorem]{Corollary}
  \newtheorem{proposition}[theorem]{Proposition}
  \newtheorem{definition}[theorem]{Definition}
\newtheorem{example}[theorem]{Example}
\newtheorem{remark}[theorem]{Remark}
\newtheorem{homework}[theorem]{Homework}
\newtheorem{case}[theorem]{}

\else
  \newtheorem{theorem}{Theorem}[section]
  \newtheorem{lemma}[theorem]{Lemma}

\fi

\usepackage[T1]{fontenc}

\journal{XXX }

\begin{document}

\begin{frontmatter}



\title{ Correlation function of a random scalar field evolving with  a rapidly fluctuating Gaussian process}

\author[2]{Jared C. Bronski}
\ead{bronski@illinois.edu}
\address[2]{Department of Mathematics, University of Illinois Urbana-Champaign, Urbana, IL 61801, United States}
\author[1]{Lingyun Ding}
\ead{dingly@live.unc.edu}
\author[1]{Richard M. McLaughlin \corref{mycorrespondingauthor}}
\cortext[mycorrespondingauthor]{Corresponding author}
\ead{rmm@email.unc.edu}
\address[1]{Department of Mathematics, University of North Carolina, Chapel Hill, NC, 27599, United States}

\begin{abstract}
  We consider a scalar field governed by an advection-diffusion equation (or a more general evolution equation) with rapidly fluctuating, Gaussian distributed random coefficients.
  In the white noise limit, we derive the closed evolution equation for the ensemble average of the random scalar field by three different strategies, i.e., Feynman-Kac formula, the limit of Ornstein-Uhlenbeck process, and evaluating the cluster expansion of the propagator on an $n$-simplex. With the evolution equation of ensemble average, we study the passive scalar transport problem with two different types of flows, a random periodic flow, and a random strain flow. For periodic flows, by utilizing the homogenization method, we show that the $N$-point correlation function of the random scalar field satisfies an effective diffusion equation at long times. For the strain flow, we explicit compute the mean of the random scalar field and show that the statistics of the random scalar field have a connection to the time integral of geometric Brownian motion.  Interestingly, all normalized moment (e.g., skewness, kurtosis) of this random scalar field diverges at long times, meaning that the scalar becomes more and more intermittent during its decay.
\end{abstract}

\end{frontmatter}

\section{Introduction}

An important and difficult problem is determining the statistics of scalar  inherits from some given external random field.
These statistics are often of great physical interest. Some examples include the understanding of the evolution of the probability distribution function for fluctuations in a passive scalar diffusing in an incompressible turbulent random flow field \cite{kraichnan1968small,sinai1989limiting,majda1993random,majda1993explicit,sukhatme2004probability,kimura1993statistics,pumir1991exponential}, the prediction of escape times for a particle
diffusing in a random fluctuating potential \cite{stein1989mean,stein1990escape}, diffusing diffusivity model \cite{chechkin2017brownian,jain2017diffusing,tyagi2017non,uneyama2019relaxation}, light propagating through random media~\cite{stephen1988temporal} and random water waves impinging on a step \cite{bolles2019anomalous,majda2019statistical} and numerous other examples.

The major question  to be addressed can be formulated as follows: Is there a universal probability distribution of a scalar in a random velocity field? If this distribution exists, then what are its properties? Naturally, the determination of these statistics requires understanding the evolution of moments of the scalar field. However, it is in general difficult to derive closed equations governing the evolution of moments since these problems involve variable coefficients.

In certain situations, however, full moment closure is available. Such situations typically involve problems whose coefficients are rapidly fluctuating in time. For the case of a passive scalar, Kraichnan \cite{kraichnan1968small,kimura1993statistics} considered a stationary isotropic turbulent velocity field admitting Gaussian white noise
time statistics and was able to derive a closed expression for the time evolution of the scalar's mean and covariance. Majda \cite{majda1993random,majda1999simplified}, studying 
the passive scalar problem in a rapidly fluctuating shear layer, obtained closed equations for the evolution of the $N$-point correlation function by
manipulating a Feynman-Kac path integral representation for the scalar's Fourier transform in the shear direction. For the special case where the shear layer is linear, these equations may be exactly solved to explicitly determine the scalar statistics inherited from the random flow and to illustrate situations leading to non-Gaussian distributions \cite{majda1993random,bronski2000problem,bronski2000rigorous,vanden2001non}.
Unfortunately, it is not at all obvious how to generalize the shear layer example to more general flow fields.

The main purpose of this paper is to demonstrate that such moment closures involving Gaussian white noise in shear layer example are special cases of a more general result. We present three strategies to derive the evolution equation of the $N$-point correlation function of the scalar field, i.e, the moment closure equation. Different strategies are suitable for different types of random flow and different governing equations. The first strategy is suitable for the advection-diffusion equation with the flow multiplied by a Gaussian white noise process. We derive the moment closure equation by manipulating the Fokker-Planck equation or Feynman-Kac formula.
The second method uses the fact that the stationary Ornstein-Uhlenbeck (OU) process converges to the Gaussian white noise process as the correlation time of the OU process vanishes.  With the moment closure equation for the OU case derived by Resnick~\cite{resnick1996dynamical}, the perturbation method yields the equation for the white noise case. The methods may be generalized to study the advection-diffusion equation with a more complicated flow that is a function of a Gaussian white noise process. 
In the third method, we consider a very general time-dependent equation with multiplicative white noise. We first derive the moment closure equation by formally using the pseudo-differential operator and the characteristic function of normal random variables. Next, we rigorously prove this result by evaluating the complete Picard iteration expansion. The assemble average of the multi-dimensional integral in the expansion is difficult to compute in general. Thanks to the correlation structure of the white noise process, we can find the rule for contracting of the multi-dimensional integral and show the expansion coverges to the moment closure equation we formally derived before.  

The paper is organized as follows. In section \ref{sec:setup}, we formulate the evolution equation of passive scalar and the evolution equation of its $N$-point correlation function. By employing the strategies we described in the previous paragraph, we derive the ensemble averaged equation of the advection-diffusion equation with rapidly fluctuating random flows in section \ref{sec:Feynman-kacFormulaBasedApproach} and \ref{sec:LimitOU}, and derive the ensemble averaged equation of a general evolution equation in section \ref{sec:PseudoOperator}. In section \ref{sec:application}, as applications of the averaged equations we derived, we explore the statistics of the random scalar field which is advected by the random strain flows and random periodic flows. In the strain flows case, we show that the statistics have a close relationship with the time integral of geometric Brownian motion which has an important application in Asian option pricing. Based on the homogenization method, we derive the effective equation of $N$-point correlation function and calculate the limiting probability distribution at long times.

\section{Preliminaries}
\label{sec:setup}
We consider the following general evolution equation with initial condition $T_{I} (\mathbf{x})$
\begin{equation}\label{eq:RandomEvolution}
\begin{aligned}
\partial_{t}  T= & A (\partial_{\mathbf{x}},t)T, \quad T (\mathbf{x},0) = & T_{I} (\mathbf{x}),\\
\end{aligned}
\end{equation}
where $A (\partial_{\mathbf{x}},t)$ is a random linear operator which depends on the derivative or higher derivative with respect to the spatial variable $\mathbf{x}$, for example, the Laplacian $\Delta_{\mathbf{x}}= \sum\limits_{j=1}^{N} \partial_{x_{j}}^{2}$.

Usually, the random scalar field $T (\mathbf{x},t)$ can be characterized by the $N$-point correlation function $\Psi_{N} (\mathbf{x}_{1},..,\mathbf{x}_{N},t) =E_{A} \left( \prod\limits_{j=1}^N T (\mathbf{x}_{j},t) \right)$ where $T (\mathbf{x}_{j},t)$ is the solution of equation \eqref{eq:RandomEvolution} with the same realization of $A$ and different initial condition $T_{I} (\mathbf{x}_{j})$. There are two reasons for that. First, the mean, variance and skewness which provides important statistical information regarding the random scalar field can be computed from the first three correlation functions.
Second, with the computed $N$-point correlation function, one can analyze the tail of PDF \cite{bronski2000problem,bronski2000rigorous}, or even obtain the analytic expression for the full PDF \cite{bronski2007explicit,ding2020ergodicity} via the conclusion of moment problem \cite{shohat1943problem}.

A nice property of equation \eqref{eq:RandomEvolution} which makes  easier to obtain the the $N$-point correlation function is that the product $\Theta (\mathbf{x}_{1},...,\mathbf{x}_{N},t) =\prod\limits_{j=1}^N T (\mathbf{x}_{j},t)$ satisfies a same type of equation with enlarged coordinates space,
\begin{equation}
\begin{aligned}
\partial_{t}  \Theta= & \sum\limits_{j=1}^{N}A (\partial_{\mathbf{x}_{j}},t)\Theta, \quad \left. \Theta \right|_{t=0 } = &\prod\limits_{j=1}^{N} T_{I} (\mathbf{x}_{j}).
\end{aligned}
\end{equation}
Hence, once we could derive the evolution equation for the mean of $T$, we could also obtain the evolution equation for the mean of $\Theta$ which is the $N$-point correlation function of the random scalar field $T$.

\section{Feynman-Kac formula based approach}
\label{sec:Feynman-kacFormulaBasedApproach}
In this section, we will study the following $d$ dimensional random advection-diffusion equation with the initial condition $T_{I}\left(\mathbf{x}\right)$,
\begin{equation}\label{eq:advectionDiffusion1}
\begin{aligned}
\partial_{t}T+\xi (t)\mathbf{v}(\mathbf{x},t)\cdot \nabla T=\kappa \Delta T, \quad  T(\mathbf{x},0)=T_{I}(\mathbf{x}),
\end{aligned}
\end{equation}
where $\kappa$ is the diffusivity, $\mathbf{v} (\mathbf{x},t)$ is a $d$-dimensional deterministic vector field, $\xi(t)$ is a  zero-mean, Gaussian random process with the correlation function given by $\left\langle \xi(t)\xi(s) \right\rangle=R(t,s)$. A special case of flow is the shear flow $\xi (t)\mathbf{v}= (u(y)\xi(t),0 )$. This type of flow can  originate from either a time varying pressure field, or by randomly moving portions of the boundary, in a high viscosity fluid \cite{ding2021enhanced}. Here, we assume that $\xi(t)$ is a Gaussian white noise in time so that $R(t,s)=g^{2} \delta(t-s)$. In most applications of this equation, the white-noise term emerges as limit of a colored noise. Therefore, we interpret this equation in the sense of Stratonovich, which will be indicated by the symbol $\circ$ in this paper.

For a fixed realization of velocity $\xi (t)\mathbf{v} (\mathbf{x},t)$, the Feynman-Kac formula yields the path integral representation for the solution of equation \eqref{eq:advectionDiffusion1},
\begin{equation}
\begin{aligned}
T(\mathbf{x},t)= & E_{ \mathbf{B} } (T_{I} (\mathbf{X}(t))), \\
\end{aligned}
\end{equation}
where the random process $\mathbf{X} (s)=\mathbf{X}^{(\mathbf{x},t)} (s)$ is the solution of stochastic differential equation (SDE),
\begin{equation}
\begin{aligned}
\mathrm{d} \mathbf{X} (s)= & - \xi (t-s)\mathbf{v} (\mathbf{X} (s),t-s)\mathrm{d}s + \sqrt{2\kappa} \mathrm{d} \mathbf{B} (s), \quad \mathbf{X} (0) = \mathbf{x}. \\
\end{aligned}
\end{equation}
Here $\mathbf{B} (t)$ is the $d$-dimensional standard Brownian motion. Notice that both the Gaussian white noise process is  stationary and temporally homogeneous, so $\xi (t-s) = \xi (s)$ in the sense of distribution. With the fact $\xi (s)\mathrm{d}  s= g\mathrm{d} \tilde{B} (s)$, where $\tilde{B} (t)$ is independent of $\mathbf{B} (t)$, we have
\begin{equation}
\begin{aligned}
\mathrm{d} \mathbf{X} (s)= & - g\mathbf{v} (\mathbf{X} (s),t-s)\circ \mathrm{d}\tilde{B} (s) + \sqrt{2\kappa} \mathrm{d} \mathbf{B} (s).\\
\end{aligned}
\end{equation}
It will be convenient to convert the Stratonovich integral into the It\^o integral,
\begin{equation}
\begin{aligned}
\mathrm{d} \mathbf{X} (s)= &g^{2}(\nabla\mathbf{v}(\mathbf{X} (s),t-s)) \mathbf{v}(\mathbf{X} (s),t-s)\mathrm{d}t - g\mathbf{v} (\mathbf{X} (s),t-s)\mathrm{d}\tilde{B} (s) + \sqrt{2\kappa} \mathrm{d} \mathbf{B} (s).\\
\end{aligned}
\end{equation}
By utilizing the Feynman-Kac formula, the ensemble average of $T (\mathbf{x},t)$ with respect to the flow, $\Psi=E_{\xi (t)}(T(x,t))= E_{\xi (t)}( E_{ \mathbf{B} } (T_{I} (\mathbf{X}(t))))$, satisfies the equation
\begin{equation}\label{eq:whiteNoiseMoment}
\begin{aligned}
 \partial_{t} \Psi=\left( \kappa \Delta+ \frac{g^{2}}{2} (\mathbf{v}\cdot \nabla)^{2} \right)    \Psi, \quad   \Psi(\mathbf{x},0)= T_{I}(\mathbf{x}).
\end{aligned}
\end{equation}
When the velocity field in equation \eqref{eq:advectionDiffusion1} is the shear flow, $\xi (t)\mathbf{v} (\mathbf{x},t) =(\xi (t) y,0)$, the advection operator  in the equation of $\Theta=\prod\limits_{j=1}^N T (\mathbf{x}_{j},t)$ is $\sum\limits_{j=1}^N y_{j} \partial_{x_{j}}$. With equation \eqref{eq:whiteNoiseMoment}, we obtain the conclusion in \cite{majda1993random}
\begin{equation}
\begin{aligned}
 \partial_{t} \Psi=\left( \kappa \Delta+ \frac{g^{2}}{2} \left(\sum\limits_{j=1}^{N}y \partial_{x_{j}}  \right)^{2} \right)    \Psi, \quad   \Psi(\mathbf{x},0)= T_{I}(\mathbf{x}).
\end{aligned}
\end{equation}

\section{Limit of Ornstein-Uhlenbeck process}
\label{sec:LimitOU}
In this section, we present the second method for deriving the moment evolution equation, which may allow us to study more general random flows. This method utilizes the fact that the stationary Ornstein-Uhlenbeck (OU) process converges to the Gaussian white noise process as the correlation time of the OU process vanishes.

We consider the advection-diffusion equation with more general random flow,
\begin{equation}\label{eq:OU1}
\begin{aligned}
\partial_{t}T+\mathbf{v}(\mathbf{x},\xi(t),t)\cdot \nabla T=\kappa \Delta T, \quad  T(x,y,0)=T_{I}(x,y). \\
\end{aligned}
\end{equation}
where $\xi(t)$ is the stationary Ornstein-Uhlenbeck process  with the damping $\gamma$ and dispersion $\sigma$ is  the solution of stochastic differential equation (SDE) $\mathrm{d}\xi (t) =-\gamma \xi (t)\mathrm{d}t +\sigma \mathrm{d}B (t)$ with initial condition $\xi(0) \sim \mathcal{N} (0, {\sigma^2}/{2 \gamma})$. Here $B (t)$ is the standard Brownian motion and $\mathcal{N} (a,b)$ is the normal distribution with mean $a$ and variance $b$.   The correlation function of $\xi (t)$ is $R(t,s)=\frac{\sigma^{2}}{2\gamma}e^{-\gamma \left| t-s \right|}$. $\gamma^{-1}$ is often referred to as the correlation time of the OU process. It is easy to check that the stationary Ornstein-Uhlenbeck process converges to the Gaussian white noise process as the correlation time vanishes with fixed $g={\sigma}/{\gamma}$.  Due to this property, we can obtain the equation for the white noise case by applying the perturbation method on the equation for the OU process case.

 Resnick \cite{resnick1996dynamical} derived the PDE for the ensemble average $\Psi=E_{\xi (t)} (T)$ by utilizing the Feynman-Kac formula. By introducing an extra variable $z$  to represent the initial value of the stationary OU process,  we have $\Psi (\mathbf{x},t)= \frac{1}{\sqrt{\pi}}\int\limits_{-\infty}^{+\infty}  \psi (\mathbf{x},z,t) e^{-z^2} \mathrm{d} z$, where $\psi (\mathbf{x},z,t)$ satisfies the following partial differential equation:
\begin{equation}\label{eq:closureeqnOU}
\begin{aligned}
\frac{\partial \psi}{\partial t}+  \mathbf{v}(\mathbf{x},g\sqrt{\gamma}z,t)\cdot \nabla  \psi + \gamma z \frac{\partial \psi}{\partial z} &= \kappa\Delta \psi+ \frac{\gamma}{2} \frac{\partial^2 \psi}{\partial z^{2}},\quad \psi(\mathbf{x},z,0)=  T_{I}(\mathbf{x}).
\end{aligned}
\end{equation}
We note that by using the same strategy, when the flow depends on $M$  independent OU processes, one can obtain the PDE for the $N$-point correlation function with $M$ new variables $z_i$, $1\leq i\leq M$. This generalization allows one to study the flow which is driven by a more general stochastic process, for example, the general stationary Gaussian random process which can be represented as a linear combination of OU processes \cite{resnick1996dynamical,csaki1991infinite}, the Cox-Ingersoll-Ross process in mathematical finance which is a sum of squared OU processes (chapter 6 of \cite{jeanblanc2009mathematical}).

For the linear shear case, $\mathbf{v}(\mathbf{x},t,\xi (t))\cdot \nabla =\xi (t)y\partial_{x}$, Resnick derived the exact expression for $\Psi$ and showed it converges to the white case in the limit  $\gamma\rightarrow \infty$ of the damping OU parameter. Here, we will generalize this conclusion. We use the velocity filed $\mathbf{v}(\mathbf{x},\xi(t),t)= \xi (t) \mathbf{v} (\mathbf{x},t)$ as an example to show the derivation of the effective equation in large $\gamma$ limit. For this type of flow, equation \eqref{eq:closureeqnOU} becomes
\begin{equation}\label{eq:closureeqnOU1}
\begin{aligned}
\frac{\partial \psi}{\partial t}+ g\sqrt{\gamma}z \mathbf{v}(\mathbf{x},t)\cdot \nabla  \psi + \gamma z \frac{\partial \psi}{\partial z} &= \kappa\Delta \psi+ \frac{\gamma}{2} \frac{\partial^2 \psi}{\partial z^{2}},\quad \psi(\mathbf{x},z,0)=  T_{I}(\mathbf{x}).
\end{aligned}
\end{equation}
Inspired by the coefficient $\sqrt{\gamma}$ in equation, we assume that the solution has the series expansion $\psi= \sum\limits_{n=0}^{\infty}\psi_n \gamma^{-\frac{n}{2}}$. Substituting the expansion into the equation \eqref{eq:closureeqnOU1} and grouping all terms of the same order of $\gamma$, we find that we have to solve the sequence of equations. The first three equations are 
\begin{equation}
\begin{aligned}
 \left( z \partial_{z}-\frac{1}{2} \partial_{z}^{2} \right)\psi_0= & 0,\\
 \left( z \partial_{z}-\frac{1}{2} \partial_{z}^{2} \right)\psi_1= &- gz \mathbf{v} (\mathbf{x},t)\cdot\nabla\psi_0, \\
 \left( z \partial_{z}-\frac{1}{2} \partial_{z}^{2} \right)\psi_2= &- gz  \mathbf{v} (\mathbf{x},t)\cdot\nabla\psi_1- \frac{\partial \psi_0}{\partial t}+\kappa \Delta \psi_{0}.\\
\end{aligned}
\end{equation}
We have $\psi_0 (\mathbf{x},t,z) =\psi_0 (\mathbf{x},t)$, $\psi_1 (\mathbf{x},t,z)=- g z \mathbf{v} (\mathbf{x},t)\cdot\nabla\psi_0 + F (\mathbf{x},t)$. The equation for $\psi_2$ becomes
\begin{equation}
\begin{aligned}
 \left( z \partial_{z}-\frac{1}{2} \partial_{z}^{2} \right)\psi_2= &g^{2}z^{2} \left( \mathbf{v} (\mathbf{x},t)\cdot\nabla \right)^{2}\psi_0-gz  \mathbf{v} (\mathbf{x},t)\cdot\nabla F- \frac{\partial \psi_0}{\partial t}+\kappa  \Delta \psi_{0} .
\end{aligned}
\end{equation}
Then $\psi_2$ takes this form 
\begin{equation}
\begin{aligned}
\psi_2= (4z^{2}-2)F_{2}+F_0- gz \mathbf{v} (\mathbf{x},t)\cdot\nabla F.
\end{aligned}
\end{equation}
Hence, we have
\begin{equation}\label{eq:OUexpansion1}
\begin{aligned}
  \psi&=\psi_0+ \frac{- gz  \mathbf{v} (\mathbf{x},t)\cdot\nabla\psi_0 + F (t,\mathbf{x})}{\sqrt{\gamma}}+ \frac{(4z^{2}-2)F_{2}+F_0- gz  \mathbf{v} (\mathbf{x},t)\cdot\nabla F}{\gamma}+\mathcal{O} (\gamma^{- \frac{3}{2}}).
\end{aligned}
\end{equation}
Since the  ensemble average $\Psi=E_{\xi (t)} (T)$ is the integral of $\psi$, $\Psi (\mathbf{x},t)= \frac{1}{\sqrt{\pi}}\int\limits_{-\infty}^{+\infty}  \psi (\mathbf{x},z,t) e^{-z^2} \mathrm{d} z$,  we take the integration at both side of equation \eqref{eq:closureeqnOU1}
\begin{equation}\label{eq:closureeqnOU2}
\begin{aligned}
\partial_t \Psi+  \frac{1}{\sqrt{\pi}}\int\limits_{-\infty}^{+\infty}  g\sqrt{\gamma}z \mathbf{v}(\mathbf{x},t)\cdot \nabla  \psi  e^{-z^2} \mathrm{d} z  &= \kappa\Delta \Psi,\quad \Psi(\mathbf{x},0)=  T_{I}(\mathbf{x}).
\end{aligned}
\end{equation}
Substituting the expansion \eqref{eq:OUexpansion1} into equation \eqref{eq:closureeqnOU2}, we have
\begin{equation}
  \begin{aligned}
    \partial_t \Psi-  \frac{g^{2}}{2}\left( \mathbf{v}(\mathbf{x},t)\cdot \nabla \right)^{2}\Psi  &= \kappa\Delta \Psi+\mathcal{O} (\gamma^{-\frac{1}{2}}),\quad \Psi(\mathbf{x},0)=  T_{I}(\mathbf{x}).
\end{aligned}
\end{equation}
By taking the limit $\gamma\rightarrow \infty$, we obtain equation \eqref{eq:whiteNoiseMoment} again.

\section{General Evolution Equation}
\label{sec:PseudoOperator}
In this section, we will show that the conclusions of advection-diffusion equations we obtained in the previous sections are special cases of a more general result. We consider  $A (\partial_{\mathbf{x}},t)= C (\partial_{\mathbf{x}},t)+ \xi (t) V (\partial_{\mathbf{x}},t)$ in equation \eqref{eq:RandomEvolution}, $C (\partial_{\mathbf{x}},t), V (\partial_{\mathbf{x}},t)$  are deterministic linear operators which depend on the derivative or higher derivative with respect to the spatial variable $\mathbf{x}$, or random linear operators which are independent of the Gaussian white noise process $\xi (t)$. Hence, equation \eqref{eq:RandomEvolution} becomes,
\begin{equation}\label{eq:GeneralRandomEvolution}
\begin{aligned}
\partial_{t}  T= & \left( C (\partial_{\mathbf{x}},t)+ \xi (t) V (\partial_{\mathbf{x}},t)\right)T, \quad T (\mathbf{x},0) = & T_{I} (\mathbf{x}).\\
\end{aligned}
\end{equation}
Obviously, equation \eqref{eq:GeneralRandomEvolution} reduces to the advection-diffusion equation \eqref{eq:advectionDiffusion1} by taking $C (\partial_{\mathbf{x}},t)=\Delta$ and $V \left(\partial_{\mathbf{x}},t)\right)=\mathbf{v} (\mathbf{x},t)\cdot \nabla $. Another interesting example which can be treated using this method is the linearized Kuramoto-Sivashinsky equation \cite{yakhot1981large}
\begin{equation}
\begin{aligned}
 &T_{t}=-T_{xx}-\epsilon T_{xxxx}+ \xi  (t) U (x,t) T_x. \\
\end{aligned}
\end{equation}
Such an equation can be analyzed in the constant coefficient case by Fourier methods \cite{avellaneda1994simple} but a new approach is needed to treat the non-constant coefficient case. Also note that path integral based approaches fail, since there is no known path integral for the above equation. 

We can formally write down the solution with the functional of operators (propagator) $U(t)$,
\begin{equation}
\begin{aligned}
T (\mathbf{x},t) = U (t) T_{I} (\mathbf{x}) =&\exp \left(\int\limits_0^t   C (\partial_{\mathbf{x}},s)+ \xi (s) V (\partial_{\mathbf{x}},s) \mathrm{d} s  \right) T_{I} \\
 =&\exp \left(\int\limits_0^t   C (\partial_{\mathbf{x}},s) \mathrm{d}s+ g\int\limits_0^{t} V (\partial_{\mathbf{x}},s) \mathrm{d} B (s)  \right) T_{I}. \\
\end{aligned}
\end{equation}
The second step follows the fact $ \xi (s) V (\partial_{\mathbf{x}},s) \mathrm{d} s =  gV (\partial_{\mathbf{x}},s) \circ \mathrm{d} B (s)= gV (\partial_{\mathbf{x}},s) \mathrm{d} B (s)$  in the sense of distribution. 
This functional calculus method was used in studying the advection-diffusion equation in \cite{taylor2012random}. Here, we formally treat the operators as ordinary variables. For a normal random variable $X \sim \mathcal{N} (\mu, \sigma)$, the characteristic function of $X$ yields $E_{X} (e^{X})=e^{\mu+\frac{1}{2}\sigma^{2}}$.  With this conclusion and It\^o isometry $E_{B} \left( \int\limits_0^{t} V (\partial_{\mathbf{x}},s) \mathrm{d} B (s) \right)^{2}= \int\limits_0^{t} V (\partial_{\mathbf{x}},s)^{2} \mathrm{d}s$, we obtain
\begin{equation}
\begin{aligned}
\Psi (\mathbf{x},t) =E_{\xi (t)} (U (t))T_{I}= &\exp \left(\int\limits_0^t  C(\partial_{\mathbf{x}},s) \mathrm{d}s+\frac{g^{2}}{2} \int\limits_0^{t} V (\partial_{\mathbf{x}},s)^{2} \mathrm{d}s   \right) T_{I}. \\
\end{aligned}
\end{equation}
From this expression, we observe that $\Psi$ satisfies the equation
\begin{equation}
\begin{aligned}
&\partial_t \Psi=  \left( C (\partial_{\mathbf{x}},t)+ \frac{g^2}{2} V (\partial_{\mathbf{x}},t)^{2}  \right)\Psi, \quad \Psi (\mathbf{x},0) =  T_{I} (\mathbf{x}).\\
\end{aligned}
\end{equation}

The above formal proof could be rigorously justified by the pseudo-differential operator. Alternatively, here, we rigorously prove this result by evaluating the Picard iteration expansion. 
The basis for our calculation is the time ordered product expansion for the propagator $U (t)$, which is given by 
\begin{equation}\label{eq:PicardRandomPropagator}
\begin{aligned}
  U (t) = &I+ \int\limits_0^t  A (\partial_{\mathbf{x}},s)\mathrm{d} s+ \int\limits_0^t \int\limits_0^{s_2}A (\partial_{\mathbf{x}},s_2)A (\partial_{\mathbf{x}},s_1)\mathrm{d}s_1\mathrm{d} s_2 \\
  &+ \int\limits_0^t \int\limits_0^{s_3}\int\limits_0^{s_2} A (\partial_{\mathbf{x}},s_3) A (\partial_{\mathbf{x}},s_2)A (\partial_{\mathbf{x}},s_1)\mathrm{d}s_1\mathrm{d} s_2\mathrm{d} s_{3}+ \hdots \\
\end{aligned}
\end{equation}
The usual practice at this point is to introduce the time ordering operator, change the integration to be over the cube of side $t$, rather than over the $n$-simplex \newline
$\Omega_{n} (t) =\left\{ (s_{1},...,s_{n})|0\leq s_1 \leq s_{2} \leq s_{3} \hdots \leq s_n \leq t \right\}$, and introduce a factor of $1/n!$. For reasons that will become clear shortly we choose to leave the series in the form written above, with the $n$-fold integration over the $n$-simplex.

Next, our goal is to show that after ensemble average with respect to $\xi (t)$, the expansion of random propagator \eqref{eq:PicardRandomPropagator} becomes
\begin{equation}\label{eq:PicardAveragedRandomPropagator}
\begin{aligned}
&E_{\xi (t)} (U (t))= I+ \int\limits_0^t  C (\partial_{\mathbf{x}},s)+\frac{g^2}{2}V (\partial_{\mathbf{x}},s)^{2}\mathrm{d} s\\
  &+ \int\limits_0^t \int\limits_0^{s_2} \left(  C (\partial_{\mathbf{x}},s_{2})+\frac{g^2}{2}V (\partial_{\mathbf{x}},s_{2})^{2} \right) \left(  C (\partial_{\mathbf{x}},s_{1})+\frac{g^2}{2}V (\partial_{\mathbf{x}},s_{1})^{2} \right)\mathrm{d} s_{1}\mathrm{d} s_{2}+... \\
\end{aligned}
\end{equation}
In doing so, we compute the ensemble average of the above propagator over the statistics of $\xi$ utilizing a cluster expansion. We use the fact that $\xi$ has Gaussian white noise statistics to ``break the bars'' in the usual way. Generally, it is difficult, if not impossible, to interpret this averaged propagator as the solution of an evolution equation. However, in the white noise limit, we will establish that this averaged propagator does in fact obey a simple evolution equation.

Since $\xi (t)$ is mean zero any terms with an odd number of factor of $V (t)$ vanish upon taking the expectation over $\xi$. The terms which appear at $n$-th order are products of $2k$ factors of $V (s_i)$ with $n-2k$ factors of $C (s_i)$, multiplied by $k$ delta functions. The arguments of the delta functions are difference between the arguments of the $V$'s. The integral is over the $n$-simplex  $\Omega_{n}$. For instance the terms appearing at fourth order are
\begin{equation}
\begin{aligned}
& C (s_4) C (s_3) C (s_2) C (s_1)+ \delta (s_4-s_3) V (s_4)V (s_3)C (s_2)C (s_1) \\
+&\delta (s_{4}-s_{2})V (s_4)C (s_3)V (s_2)C (s_1)+\delta (s_{4}-s_{1}) V (s_4)C (s_3)C (s_2)V (s_1) \\
+&\delta (s_{3}-s_{2})C (s_4)V (s_3)V (s_2)C (s_1)+\delta (s_{3}-s_{1}) C (s_4)V (s_3)C (s_2)V (s_1) \\
+&\delta (s_{2}-s_{1})C (s_4)C (s_3)V (s_2)V (s_1)\\
+& \left(\delta (s_{4}-s_{3})\delta (s_{2}-s_{1})
  +\delta (s_{4}-s_{2})\delta (s_{3}-s_{1})+\delta (s_{4}-s_{1})\delta (s_{3}-s_{2})  \right) V (s_4)V (s_3)V (s_2)V (s_1). \\
\end{aligned}
\end{equation}

The utility of leaving the integration over the simplicity regions now becomes apparent. The integral over the delta functions can be done easily. The result of this integration is, of course, the characteristic function of the support of the delta function. The important thing to realize for the purposes of this calculation is that the intersection of the $n-1$ dimensional hyperplane $s_i=s_{i+k}$ with the $n$-simplex  $0\leq s_1 \leq s_{2} \leq s_{3} \hdots \leq s_n \leq t$ is the $n-k$ dimensional simplex $0\leq s_1 \leq ... \leq s_{i}=s_{i+1}=...=s_{i+k} \leq  \hdots \leq s_n \leq t$. Thus the intersection of the support of the delta function with the region of integration is a region of dimension $n-k$. Integrating the characteristic function of an $n-k$ dimensional set over an $n-1$ dimensional set gives no contribution unless $k=1$. The upshot is that the only contribution which do not vanish are those from adjacent factors of $V$ in the product. It is precisely the non-adjacent terms which make it difficult to interpret the averaged propagator as the solution to an evolution equation. In the white noise limit, these terms do not contribute. In addition, we should also notice that since the support of delta function lies on the boundary of $\Omega_n (t)$, the integration produces the prefactor $1/2$. We summarize the above discussion as the following formula
\begin{equation}\label{eq:IntDeltaOnSimplex}
\begin{aligned}
\int\limits_{\Omega_n (t)}^{}\delta (s_{i}-s_j)f (\mathbf{s}_{n})\mathrm{d} \mathbf{s}_{n}= &
\begin{cases}
  \frac{1}{2}\int\limits_{\Omega_{n-1} (t)}^{}f (s_{1},...,s_{j-1},s_{i},s_{j+1},...,s_{j})\mathrm{d} \mathbf{s}_{n-1}, & i=j+1,\\
  0,&\text{otherwise,}\\
\end{cases}
\end{aligned}
\end{equation}
where $i>j$, $\mathbf{s}_n= (s_{1},s_{2},...,s_{n})$. We also relabel $s_n$ as $s_j$ in the right hand side integral to make the formula more concise. 

It is worthwhile to verify the above claim for the first several terms. After taking ensemble average, the first time integral in \eqref{eq:PicardRandomPropagator} becomes
\begin{equation}
\begin{aligned}
E_{\xi(t)}\left( \int\limits_0^t  A (\partial_{\mathbf{x}},s)\mathrm{d} s \right) =& \int\limits_0^t C (\partial_{\mathbf{x}},s)\mathrm{d} s. \\
\end{aligned}
\end{equation}
The second term becomes
\begin{equation}
\begin{aligned}
  &E_{\xi(t)}\left(\int\limits_0^t \int\limits_0^{s_2}A (\partial_{\mathbf{x}},s_2)A (\partial_{\mathbf{x}},s_1)\mathrm{d}s_1\mathrm{d} s_2  \right) = \int\limits_0^t \int\limits_0^{s_2} C (\partial_{\mathbf{x}},s_{2})C (\partial_{\mathbf{x}},s_{1})\mathrm{d} s_{1}\mathrm{d} s_{2} +\int\limits_0^t \frac{g^2}{2}V (\partial_{\mathbf{x}},s)^{2}\mathrm{d} s.\\
\end{aligned}
\end{equation}
The third term becomes
\begin{equation}
  \begin{aligned}
    &E_{\xi(t)}\left(\int\limits_{\Omega_3(t)}^{}\prod\limits_{i=1}^{3} A (\partial_{\mathbf{x}},s_i)\mathrm{d} \mathbf{s}_{3}\right) =\int\limits_{\Omega_3(t)}^{}\prod\limits_{i=1}^{3} C (\partial_{\mathbf{x}},s_i)\mathrm{d} \mathbf{s}_{3}+  \frac{g^2}{2}\int\limits_0^t \int\limits_0^{s_2} C (\partial_{\mathbf{x}},s_{2})V (\partial_{\mathbf{x}},s_{1})^{2}\mathrm{d} s_{1}\mathrm{d} s_{2} \\
  &+  \frac{g^2}{2}\int\limits_0^t \int\limits_0^{s_2} V (\partial_{\mathbf{x}},s_{2})^{2}C (\partial_{\mathbf{x}},s_{1})\mathrm{d} s_{1}\mathrm{d} s_{2}. \\
\end{aligned}
\end{equation}

Now, with these examples and conclusion\eqref{eq:IntDeltaOnSimplex}, we can show that equation \eqref{eq:PicardRandomPropagator} becomes equation \eqref{eq:PicardAveragedRandomPropagator} after taking ensemble average. After taking ensemble average,  any term in the cluster expansion of equation \eqref{eq:PicardRandomPropagator} is a a integral of the product of $C (s_i)$ and $\frac{g^{2}}{2}V (s_j)$. Clearly, that is a term in \eqref{eq:PicardAveragedRandomPropagator}. On the other hand, for any term in  \eqref{eq:PicardAveragedRandomPropagator} we can find the unique corresponding term in equation \eqref{eq:PicardRandomPropagator} by writing $\frac{g^2}{2}V (s_{i})^{2}$ as $\xi (s_{i}) V (s_{i})\xi (s_{i+1}) V (s_{i+1})$ and suitable relabeling the subscripts.

Last, it should be stressed that the crucial element in this Calculation is the fact that the bracket of non-adjacent $V (s_{i})$'s gives no contribution after integration. This observation allows us to interpret the averaged series as the propagator for the above problem. In a more general situation, for
instance one where $\xi(t)$ is Gaussian with some finite time correlation, terms like
\begin{equation}
\begin{aligned}
\int\limits_{0}^{t} \int\limits_0^{s_3} \int\limits_{0}^{s_{2}}V (s_3) C (s_2) V (s_1) \mathrm{d} s_1 \mathrm{d} s_2 \mathrm{d} s_3, \\
\end{aligned}
\end{equation}
have a non-vanishing contribution. Such terrns do
not appear in the series for the propagator in equation \eqref{eq:PicardRandomPropagator}, and prevent us from interpreting the resulting series as a time-ordered exponential.

\section{Application}
\label{sec:application}
In this section, we  apply our general result to the  case of a passive scalar diffusing in a rapidly fluctuating fluid  flow and study closed equations for the scalar's $N$-point correlation function. We are interested in two types of the flow, the periodic flow and strain flow. 

\subsection{Random periodic flow}
For the shear flow $\xi (t) \mathbf{v} (\mathbf{x},t)= \xi (t) (u (y),0)$ and the parallel-plate channel domain $\Omega=\left\{ (x,y)| x\in \mathbb{R}, y\in [0,1] \right\}$ with no-flux boundary conditions,  Camassa et al.\cite{camassa2021persisting} proved that the $N$-point correlation function satisfies an effective heat equation at long times. Interestedly, in this example, the limiting distribution of the scalar field is not Gaussian and exhibit intermittency. In fact, it is negatively skewed at sufficiently low P\'eclet number and could be positively skewed at large P\'eclet number.

For the parallel-plate channel domain, one can convert the no-flux boundary condition into periodic boundary condition by evenly extending the scalar field. Hence, we expect the similar result holds for a general periodic velocity field. We assume $\mathbf{v} (\mathbf{x},t)$ is  a periodic incompressible velocity field  with the spatial period $L_x$ and temporal period $L_{t}$.
We consider equation \eqref{eq:advectionDiffusion1} with the slow varying initial condition $T (\mathbf{x},0)= T_{I} (\frac{\mathbf{x}}{a})$, $a \gg 1$. Notice that $g^{2} \sim \text{Time}^{-1}$. With the change of variables, 
\begin{equation}
\begin{aligned}
&a\mathbf{x}'= \mathbf{x}, \quad \frac{a^2}{\kappa}t'=t,  \quad U= Lg^{2}, \quad \epsilon=\frac{L_{x}}{a},\\
& U\xi' (t') \mathbf{v}'(\mathbf{x}',t')=\xi (t) \mathbf{v} (\mathbf{x}, t), \quad \frac{g \sqrt{\kappa}}{a}\xi'(t')=\xi (t), \\
\end{aligned}
\end{equation}
we can drop the primes without confusion and obtain the nondimensionalized version of \eqref{eq:advectionDiffusion1}:
\begin{equation}\label{eq:advectionDiffusionNonDimension}
 \partial_{t}T+\text{Pe} \xi (t) \mathbf{v}\left( \frac{\mathbf{x}}{\epsilon},\frac{t}{\epsilon^{2}} \right)\cdot \nabla T=\Delta T,\quad T (\mathbf{x},0)= T_{I} (\mathbf{x}),
\end{equation}
where $\mathbf{v} (\mathbf{x},t)$ has the spatial period $1$ and temporal period $\kappa L_{t}/L_{x}^{2}$, and we have introduced the P\'{e}clet number  $\text{Pe}=  {U L}/{ \kappa}=  {L^{2}g^{2}}/\kappa$. With equation \eqref{eq:whiteNoiseMoment}, the ensemble average of the scalar field $\Psi=E_{\xi (t)} (T)$ satisfies the equation
\begin{equation}
  \begin{aligned}
     \partial_{t}\Psi= \left( \Delta + \mathrm{Pe}^{2}\left( \mathbf{v}\left( \frac{\mathbf{x}}{\epsilon},\frac{t}{\epsilon^{2}} \right)\cdot \nabla \right)^{2} \right)\Psi,\quad \Psi (\mathbf{x},0)= T_{I} (\mathbf{x}).
\end{aligned}
\end{equation}

We follow the similar multiscale analysis procedure described in the appendix 7.1 of \cite{ding2021enhanced} or \cite{mclaughlin1994turbulent}. We seek the asymptotic approximation to $\Psi(\mathbf{x},t)$ in the limit $\epsilon\rightarrow 0$ that has the following multiscale expansion,
\begin{equation}\label{eq:multiscaleExpansion}
\Psi(\mathbf{x},t)=\Psi_{0}(\mathbf{x},\mathbf{y},t,\tau)+\epsilon\Psi_{1}(\mathbf{x},\mathbf{y},t,\tau)+\epsilon^{2}\Psi_{2}(\mathbf{x},\mathbf{y},t,\tau)+\mathcal{O}(\epsilon^3),
\end{equation}
with two different scales in $\mathbf{x}$-direction: $\mathbf{x}$ (slow), $\mathbf{y}= {\mathbf{x}}/{\epsilon}$ (fast), and in $t$-direction: $t$ (slow), $\tau= {t}/{\epsilon^2}$ (fast), respectively. We use the angle bracket to denote the average of a function,
\begin{equation}
\begin{aligned}
\left\langle \mathbf{v} (\mathbf{x},t) \right\rangle_{\mathbf{x},t}= \int\limits_0^{L_{t}}\int\limits_{[0,L_{x}]^{d}}^{} \mathbf{v} (\mathbf{x},t)\mathrm{d} \mathbf{x} \mathrm{d} t.
\end{aligned}
\end{equation}
By taking the average on both side of equation \eqref{eq:multiscaleExpansion}, we have 
\begin{equation}
\begin{aligned}
\left\langle \Psi(\mathbf{x},t) \right\rangle_{\mathbf{y},\tau}=\left\langle  \Psi_{0}(\mathbf{x},\mathbf{y},t,\tau)\right\rangle_{\mathbf{y},\tau}+\epsilon \left\langle \Psi_{1}(\mathbf{x},\mathbf{y},t,\tau) \right\rangle_{\mathbf{y},\tau}+\epsilon^{2} \left\langle \Psi_{2}(\mathbf{x},\mathbf{y},t,\tau) \right\rangle_{\mathbf{y},\tau}+\mathcal{O}(\epsilon^3).
\end{aligned}
\end{equation}
Since equation \eqref{eq:multiscaleExpansion} holds for arbitrarily small $\epsilon$, we have $\left\langle \Psi(\mathbf{x},t) \right\rangle_{\mathbf{y},\tau}=\left\langle  \Psi_{0}(\mathbf{x},\mathbf{y},t,\tau)\right\rangle_{\mathbf{y},\tau}$, and  $\left\langle  \Psi_{n}(\mathbf{x},\mathbf{y},t,\tau)\right\rangle_{\mathbf{y},\tau}=0$, $n\geq 1$. In addition, the differential operators along the $\mathbf{x}$ and $t$ direction will be replaced,
\begin{equation}
\begin{aligned}
& \nabla_{\mathbf{x}}\rightarrow   \nabla_{\mathbf{x}}+\frac{1}{\epsilon} \nabla_{\mathbf{y}}, \quad  \Delta_{\mathbf{x}}\rightarrow \Delta_{\mathbf{x}}+\frac{2}{\epsilon} \nabla_{\mathbf{x}}\cdot\nabla_{\mathbf{y}}+ \frac{1}{\epsilon^2} \Delta_{\mathbf{y}} , \quad \partial_{t}\rightarrow   \partial_{t}+\frac{1}{\epsilon^{2}} \partial_{\tau}.   
\end{aligned}
\end{equation}
We would have a hierarchy of equations such that the expansion of $\Psi (\mathbf{x},t)$ holds for arbitrarily small $\epsilon$. 
For $\mathcal{O}(\epsilon^{-2})$,  we have:
\begin{equation}
\mathcal{L}\Psi_{0}= 0, \quad
 \Psi_0(\mathbf{x},\mathbf{y},t,\tau)|_{t=0,\tau=0}= T_{I}(\mathbf{x}), 
\end{equation}
where  $\mathcal{L}\Psi=\left(\partial_{\tau}- \Delta_{\mathbf{y}}- \mathrm{Pe}^{2} \left( \mathbf{v}(\mathbf{y},\tau)\cdot\nabla_{\mathbf{y}} \right)^{2} \right)\Psi$. Since the initial condition is a function of the variable $\mathbf{x}$ only, we have $ \Psi_0(\mathbf{x},\mathbf{y},t,\tau)=\Psi_0(\mathbf{x},t)
$.

For $\mathcal{O}(\epsilon^{-1})$, we have
\begin{equation}
\begin{aligned}
  &\mathcal{L}\Psi_1= 2\nabla_{\mathbf{x}}\cdot \nabla_{\mathbf{y}} \Psi_{0}+\mathrm{Pe}^{2} (\mathbf{v}(\mathbf{y},\tau)\cdot \nabla_{\mathbf{x}}) ( \mathbf{v}(\mathbf{y},\tau)\cdot \nabla_{\mathbf{y}})\Psi_{0}+\mathrm{Pe}^{2} (\mathbf{v}(\mathbf{y},\tau)\cdot \nabla_{\mathbf{y}}) ( \mathbf{v}(\mathbf{y},\tau)\cdot \nabla_{\mathbf{x}})\Psi_{0} ,\\
  &  \Psi_1(\mathbf{x},\mathbf{y},0,0)= 0. \\
\end{aligned}
\end{equation}
The first and second term on the right hand side are zero. Fredholm solvability states that the equation $\mathcal{L}\Psi=f$ has a solution if and only if $\left\langle fg \right\rangle=0$ for any solution of equation $\mathcal{L}^{*}g=0$, where $\mathcal{L}^{*}$ is the adjoint operator of $\mathcal{L}$. Here, the solvability condition is guaranteed by the periodicity of $\mathbf{v}$, $\left\langle (\mathbf{v} (\mathbf{y},\tau)\cdot \nabla_{\mathbf{y}})\mathbf{v}(\mathbf{y},\tau)  \right\rangle_{\mathbf{y},\tau}=\frac{1}{2}\left\langle \nabla_{\mathbf{y}} (\mathbf{v}(\mathbf{y},\tau)\cdot \mathbf{v} (\mathbf{y},\tau))  \right\rangle_{\mathbf{y},\tau}  =\mathbf{0}$. Due to the linearity of the equation, the general form of the solution is $\Psi_1= \boldsymbol{\theta}(\mathbf{y},\tau)\cdot \nabla_{\mathbf{x}}\Psi_{0}+C(\mathbf{x},t)$. Therefore, we have
\begin{equation}\label{eq:cellproblem}
\mathcal{L}\theta_{i}=  \mathrm{Pe}^{2}(\mathbf{v}(\mathbf{y},\tau)\cdot \nabla_{\mathbf{y}})  v_{i}(\mathbf{y},\tau),\quad \theta_{i}(\mathbf{y},0)=0, \quad 1\leq i\leq d,
\end{equation}
where $\theta_i$, $v_i$ are $i$th component of $\boldsymbol{\theta}$, $\mathbf{v}$ respectively. 
In addition, we have an extra condition $\left\langle \boldsymbol {\theta}(\mathbf{y},\tau)\right\rangle_{\mathbf{y},\tau}=\mathbf{0}$ from $\left\langle  \Psi_{1}(\mathbf{x},\mathbf{y},t,\tau)\right\rangle_{\mathbf{y},\tau}=0$. For $\mathcal{O}(\epsilon^0)$, we have
\begin{equation}
\begin{aligned}
  &  \mathcal{L}\Psi_2=  - \partial_{t}\Psi_{0}+ \left( \Delta_{\mathbf{x}} + \mathrm{Pe}^{2}\left( \mathbf{v}(\mathbf{y},\tau) \cdot \nabla_{\mathbf{x}} \right)^{2} \right)\Psi_{0}+2\nabla_{\mathbf{x}}\cdot \nabla_{\mathbf{y}} \Psi_{1}\\
  &\hspace{1.2cm}+\mathrm{Pe}^{2} (\mathbf{v}(\mathbf{y},\tau)\cdot \nabla_{\mathbf{x}}) ( \mathbf{v}(\mathbf{y},\tau)\cdot \nabla_{\mathbf{y}})\Psi_{1}+\mathrm{Pe}^{2} (\mathbf{v}(\mathbf{y},\tau)\cdot \nabla_{\mathbf{y}}) ( \mathbf{v}(\mathbf{y},\tau)\cdot \nabla_{\mathbf{x}})\Psi_{1},\\
  &\Psi_2(\mathbf{x},\mathbf{y},0,0)= 0. \\
\end{aligned}
\end{equation}
Due to the incompressibility and periodicity of the velocity field, the third term and last term on the right hand side of the above equation is zero under the average.  The solvability condition yields the effective diffusion equation with $\epsilon\rightarrow 0$
\begin{equation}\label{eq:effdiffusivity}
\begin{aligned}
  \partial_{t}\Psi_{0}&=\nabla_{\mathbf{x}} \cdot \left(\Lambda_{d}  \nabla_{\mathbf{x}}\Psi_{0}  \right), \quad (\Lambda_{d})_{i,j}= \delta_{i,j}+\mathrm{Pe} \left\langle v_{i}v_{j}+ v_{i}\sum\limits_{k=1}^{d}v_{k}\partial_{y_{k}}\theta_{j} \right\rangle_{\mathbf{y},\tau},\\
  \Lambda_{d}&=\mathbf{I}_{d}+\mathrm{Pe}^{2} \left\langle \mathbf{v}(\mathbf{y}, \tau)^{T}\mathbf{v}(\mathbf{y}, \tau)+\mathbf{v}(\mathbf{y}, \tau)^{T}  \left( \mathbf{v}(\mathbf{y}, \tau) \nabla_{\mathbf{y}} \right)\boldsymbol {\theta}(\mathbf{y},\tau)\right\rangle_{\mathbf{y},\tau}, 
\end{aligned}
\end{equation}
where $\mathbf{I}_d$ is the $d\times d$  identity matrix. $\delta_{i,j}$ is  Kronecker delta, $\delta_{i,j}=1$ if $i=j$, otherwise, $\delta_{i,j}=0$. Due to molecular diffusion, the solution of equation \eqref{eq:advectionDiffusion1} with a non-slowly varying initial condition will be a slowing varying function eventually. Hence, equation \eqref{eq:effdiffusivity} is also the effective equation of equation \eqref{eq:advectionDiffusion1} with a general initial condition at the diffusion time scale.

Next, we will compute the $N$-point correlation function for a 2 dimensional scalar field. In doing so, as we discussed in section \ref{sec:setup}, we need to consider the problem in higher spatial dimensions.
We assume the two dimensional flow is $(u_{1} (x,y,t),v_{1} (x,y,t))$. The $N$-point correlator will live in the space  $\mathbf{z}=(x_{1},y_{1},x_{2},y_{2},...,x_{N},y_{N})\in \mathbb{R}^{2N}$.  The flow is
\begin{equation}
\begin{aligned}
\mathbf{v} (\mathbf{z},t) = (u_{1} (x_{1},y_{1},t),v_{1}(x_{1},y_{1},t),u_{2} (x_{2},y_{2},t),v_{2}(x_{2},y_{2},t),...,u_{N} (x_{N},y_{N},t),v_{N}(x_{N},y_{N},t)). \\
\end{aligned}
\end{equation}
Then, we can apply the results we derived in previous section to this flow directly. 
The cell problem \eqref{eq:cellproblem} becomes
\begin{equation}\label{eq:cellproblem2D}
\begin{aligned}
  &\left(\partial_{t}-\partial_{x_{1}}^{2}-\partial_{y_{1}}^{2}-\mathrm{Pe} (u \partial_{x_{1}}+v\partial_{y_{1}})^{2}  \right)\theta_{1}=  \mathrm{Pe}^{2}(u \partial_{x_{1}}+v\partial_{y_{1}} )  u 
  ,\quad \theta_{1}(\mathbf{z},0)=0,\\
  &\left(\partial_{t}-\partial_{x_{1}}^{2}-\partial_{y_{1}}^{2}-\mathrm{Pe} (u \partial_{x_{1}}+v\partial_{y_{1}})^{2}  \right)\theta_{2}=  \mathrm{Pe}^{2}(u \partial_{x_{1}}+v\partial_{y_{1}} )  v, \quad\theta_{2}(\mathbf{z},0)=0.\\
\end{aligned}
\end{equation}
Notice that $\theta_{1},\theta_{3},...\theta_{2N-1}$ are identical under the permutation of variables. More precisely, we have
\begin{equation}
\begin{aligned}
  \theta_{3} (x_{1},y_{1},x_{2},y_{2},...,x_{N},y_{N})= &\theta_{1} (x_{2},y_{2},x_{1},y_{1},...,x_{N},y_{N}),\\
  ...&\\
  \theta_{2N-1} (x_{1},y_{1},x_{2},y_{2},...,x_{N},y_{N})= &\theta_{1} (x_{N},y_{N},x_{2},y_{2},...,x_{1},y_{1}).
\end{aligned}
\end{equation}
This properties also holds for  $\theta_{2},\theta_{4},...,\theta_{2N}$. With this fact, once we solve the cell problem \eqref{eq:cellproblem2D} for $\theta_{1}, \theta_{2}$, we can determine the leading order approximation of  $N$-point correlation function for arbitrary $N$. The diffusion tensor in the effective equation \eqref{eq:effdiffusivity} now becomes a $N\times N$ block matrix with $N^{2}$ $2\times 2$ matrix as blocks. More accurately, when $i\neq j$, we have
\begin{equation}
\begin{aligned}
\left( \Lambda_{2N} \right)_{2i-1,2j-1}=&\mathrm{Pe}^{2} \left( \left\langle u \right\rangle_{\mathbf{y},\tau}^{2}+\left\langle u \right\rangle_{\mathbf{y},\tau} \left\langle (u\partial_{x}+v\partial_{y})\theta_{1} \right\rangle_{\mathbf{y},\tau} \right),\\
\left( \Lambda_{2N} \right)_{2i-1,2j}=&\mathrm{Pe}^{2} \left( \left\langle u \right\rangle_{\mathbf{y},\tau} \left\langle v \right\rangle_{\mathbf{y},\tau}+\left\langle u \right\rangle_{\mathbf{y},\tau} \left\langle (u\partial_{x}+v\partial_{y})\theta_{2} \right\rangle_{\mathbf{y},\tau} \right),\\
\left( \Lambda_{2N} \right)_{2i,2j-1}=&\mathrm{Pe}^{2} \left( \left\langle v \right\rangle_{\mathbf{y},\tau}\left\langle u \right\rangle_{\mathbf{y},\tau}+\left\langle v \right\rangle_{\mathbf{y},\tau} \left\langle (u\partial_{x}+v\partial_{y})\theta_{1} \right\rangle_{\mathbf{y},\tau} \right),\\
\left( \Lambda_{2N} \right)_{2i,2j}=&\mathrm{Pe}^{2} \left( \left\langle v \right\rangle_{\mathbf{y},\tau}^{2}+\left\langle v \right\rangle_{\mathbf{y},\tau} \left\langle (u\partial_{x}+v\partial_{y})\theta_{2} \right\rangle_{\mathbf{y},\tau} \right).\\
\end{aligned}
\end{equation}
When $i=j$, we have
\begin{equation}
\begin{aligned}
\left( \Lambda_{2N} \right)_{2i-1,2i-1}=&1+\mathrm{Pe}^{2} \left( \left\langle u^{2} \right\rangle_{\mathbf{y},\tau}+\left\langle u(u\partial_{x}+v\partial_{y})\theta_{1} \right\rangle_{\mathbf{y},\tau} \right),\\
\left( \Lambda_{2N} \right)_{2i-1,2i}=&\mathrm{Pe}^{2} \left( \left\langle uv \right\rangle_{\mathbf{y},\tau} +\left\langle u(u\partial_{x}+v\partial_{y})\theta_{2} \right\rangle_{\mathbf{y},\tau} \right),\\
\left( \Lambda_{2N} \right)_{2i,2i-1}=&\mathrm{Pe}^{2} \left( \left\langle vu \right\rangle_{\mathbf{y},\tau}+\left\langle v(u\partial_{x}+v\partial_{y})\theta_{1} \right\rangle_{\mathbf{y},\tau} \right),\\
\left( \Lambda_{2N} \right)_{2i,2i}=&1+\mathrm{Pe}^{2} \left( \left\langle v^{2} \right\rangle_{\mathbf{y},\tau}+ \left\langle v(u\partial_{x}+v\partial_{y})\theta_{2} \right\rangle_{\mathbf{y},\tau} \right).\\
\end{aligned}
\end{equation}

When the two dimensional flow is a shear flow $(u_{1} (x,y,t),v_{1} (x,y,t))= (u (y,t),0)$, the cell problem \eqref{eq:cellproblem2D} only has a trivial solution  $\theta_{1}=\theta_{2}=0$. The diffusion tensor would be extremely simplified.  In this case, we have $  \Lambda_{2N}=\mathbf{I}_{2N}+\mathrm{Pe}^{2} \left\langle \mathbf{v}(\mathbf{y}, \tau)^{T}\mathbf{v}(\mathbf{y}, \tau)  \right\rangle_{\mathbf{y},\tau}$, which reproduces the result in \cite{camassa2021persisting,ding2020ergodicity,bronski1997scalar}. Then one can compute the PDF of the scalar field by using the procedure described in \cite{ding2020ergodicity,bronski2007explicit}. We remark that \cite{camassa2021persisting,ding2020ergodicity} only consider the case where $u (y)$ is independent of $t$. Our result presented in this section generalizes it to the time dependent case, $u (y,t)$.

\subsection{Random strain flow and time integral of geometric Brownian motion}
It is well known that any velocity field can be locally decomposed into the sum of a shear flow, a strain flow and a rotational flow. Many studies focus on the random shear flow  \cite{majda1993random,mclaughlin1996explicit}, few studies have addressed the strain flow or rotational flow. Here we apply our general results to the 2-dimensional random strain flow which takes the form $\xi (t)\mathbf{v} (\mathbf{x},t)=\xi (t) ( x,- y)$.

The mean of scalar field, $\Psi_{1} (x,y,t) =E_{\xi (t)}[(T (x,y,t)]$, satisfies the equation
\begin{equation}\label{eq:StrainMoment}
\begin{aligned}
 \partial_t \Psi_{1}=\left( \kappa \Delta+ \frac{g^{2}}{2} \left( x \partial_{x}-y\partial_{y} \right)^{2} \right)    \Psi_{1}, \quad   \Psi_{1}(x,y,0)= T_{I}(x,y).
\end{aligned}
\end{equation}
We can derive the exact solution for a special initial condition, $T_I (x,y) =\delta (x)$. In this case, the solution is independent of $y$, the equation \eqref{eq:StrainMoment} becomes
\begin{equation}\label{eq:StrainMomentN1}
\begin{aligned}
 \partial_t \Psi_{1}=\left( \kappa\Delta+ \frac{g^{2}}{2} ( x\partial_{x})^{2} \right)    \Psi_{1}, \quad   \Psi_{1}(x,0)= \delta (x).
\end{aligned}
\end{equation}
Interestingly, when $\kappa=0$, $f=e^{\frac{1}{2}g^{2}t}\Psi_1 (x,-t)$ is the solution of Black-Scholes equation,
\begin{equation}
\begin{aligned}
&\partial_{t} f+ \frac{1}{2}\sigma^{2}x^2 \partial_{x}^{2} f+ r x \partial_{x}f-rV=0, \\
\end{aligned}
\end{equation}
with the volatility $\sigma=g$ and the risk-free interest rate $r=\frac{1}{2}g^{2}$. Shortly, we will see that the random strain flow has a more close relationship to the mathematical finance. 

With the change of variable $z= \int\limits_0^x \frac{\mathrm{d} s}{\sqrt{\frac{g^{2}s^{2}}{2}+\kappa } }=\frac{\sqrt{2}}{g}\sinh ^{-1}\left(\frac{ gx}{\sqrt{ 2\kappa }}\right)$, equation \eqref{eq:StrainMomentN1} becomes the heat equation
\begin{equation}
\begin{aligned}
 \partial_{t}\Psi_{1} =\partial_{z}^{2} \Psi_{1}, \quad \Psi_{1}(z,0)= \delta (z)/\sqrt{\kappa}.
\end{aligned}
\end{equation}
Hence the solution of equation \eqref{eq:StrainMomentN1} is
\begin{equation}\label{eq:StrainMomentN1Sol}
\begin{aligned}
\Psi_{1} (x,t)=\frac{1}{2 \sqrt{\pi \kappa t} } \exp \left( -\frac{\left( \frac{\sqrt{2}}{g}\sinh ^{-1}\left(\frac{ gx}{\sqrt{ 2\kappa }}\right) \right)^2}{4 t} \right).
\end{aligned}
\end{equation}

We can study this random strain flow problem via an alternative approach. With a fixed realization of white noise $\xi (t)$, we can solve equation \eqref{eq:advectionDiffusion1} with the strain flow $\mathbf{v}=(x, -y)$ and the line source initial data $T_{I} (x)=\delta (x-x_{0})$ by utilizing the Fourier transform and  method of characteristics,
\begin{equation}
\begin{aligned}
T (x,t)= & \frac{1}{ \sqrt{4 \kappa \pi \int_0^t e^{-2 g B(s)}
           \mathrm{d} s}} \exp \left(-\frac{e^{-2 g B(t)} \left(x-x_{0}
           e^{g B(t)}\right)^2}{4 \kappa \int_0^t e^{-2 gB(s)} \mathrm{d}
           s}\right).\\
\end{aligned}
\end{equation}
First, we consider a special case of this random field. When $x=0, x_{0}=0$, we have $X_{t}=T (0,t)=\frac{1}{ \sqrt{ 4\kappa \pi \int_0^t e^{-2 g B(s)}   \mathrm{d} s}}$. This random variable is related to the time integral of geometric Brownian motion. In mathematical finance, the stock price is usually modeled by the geometric Brownian motion.  For Asian options, the payoff is determined by the average underlying price over some pre-set period of time. The time integral of geometric Brownian motion is directly related to the time average of the stock price which makes it important in Asian option pricing \cite{matsumoto2005exponential}.

Some properties of the time integral of geometric Brownian motion could help us understand the statistical properties of scalar field. Let's note $B_{t}^{(\mu)}=B_t+\mu t$, $A_{t}^{(\mu)}= \int_0^t e^{2 B_s^{\mu}}\mathrm{d} s $. We omit the superscript when $\mu=0$.  With this notation,  the random variable we interested in is $X_{t}= g\left( 4 \pi A_{g^{2}t} \right)^{-\frac{1}{2}}$. Using the time reversal and It\^o's formula (or otherwise), we have the recurrence relation between negative moments of $A_{t}^{\mu}$ \cite{donati2000positive,dufresne2001integral}.
\begin{lemma}\label{thm:Matsumoto1}
If $2m-\mu\geq 0$, it holds that
  \begin{equation}
 E[(A^{(\mu)}_t)^{-(m+1)}]= 2(m-\mu)E[(A_t^{(\mu)})^{-m}]- \frac{1}{m} \frac{\partial }{\partial t} E[(A_t^{(\mu)})^{-m}].
\end{equation}
\end{lemma}
Dufresne \cite{dufresne2000laguerre} obtained the following formula of the moments of $A_{t}^{\mu}$,
       \begin{lemma}\label{thm:dufresne2000P}
         For any $t>0$ and $r \in \mathbb{C}$ with $\Re (r)> - \frac{3}{2}$, it holds that 
         \begin{equation}
           E[(2A^{(\mu)}_t)^{r}]= \frac{e^{-\frac{\mu^2 t}{2}}}{\sqrt{2 \pi t^3}} \int\limits_0^{\infty} y e^{- \frac{y^2}{2t}} \varphi_{\mu} (r,y)\mathrm{d} y,
         \end{equation}
where, for the Gauss hypergeometric function $_{2}F_{1}$ and 
\begin{equation}
\begin{aligned}
&\varphi_{\mu} (r,y)= \frac{\Gamma (r+1)}{\Gamma (2r+2)} e^{-\mu y} (1- e^{-2y})^{1+2r} {}_{2}F_{1} (\mu+2r+1,1+r; 2+2r;1-e^{-2y}),\\
&\varphi_{\mu} (-1,y)= \frac{\cosh ((\mu-1)y)}{\sinh (y)}.
\end{aligned}
\end{equation}
       \end{lemma}
Donati-Martin et al. \cite{donati2000positive} proved that  
\begin{lemma}\label{thm:Matsumoto}
  For any $t>0$, it holds that
  \begin{equation}
\begin{aligned}
E \left[ A_{t}^{-1} \right]= & \int\limits_0^{\infty} k e^{-k^{2}t/2} \coth \left( \frac{\pi}{2}k \right)\mathrm{d} k.\\
\end{aligned}
\end{equation}
\end{lemma}

Using lemma \ref{thm:dufresne2000P} or Bougerol's identity \cite{bertoin2013some}, we have $E[A_{t}^{-\frac{1}{2}}]= t^{-\frac{1}{2}}$. Hence, $E[X_t]= g\left( 4 \pi \right)^{-\frac{1}{2}} \left( g^{2}t \right)^{-\frac{1}{2}}= (4\pi t)^{-\frac{1}{2}}$, which is consistent with equation \eqref{eq:StrainMomentN1Sol}. With lemma \ref{thm:Matsumoto} and Laplace's method \cite{kirwin2010higher}, we obtain the asymptotic expansion of $E [A_t^{-1}]$ as $t\rightarrow \infty$,
\begin{equation}
\begin{aligned}
E[A_{t}^{-1}] &=\sqrt{\frac{2}{\pi t}}+  \frac{\pi ^{3/2}}{6 \sqrt{2} t^{3/2}}-\frac{\pi ^{7/2}}{120 \sqrt{2} t^{5/2}}+ \mathcal{O} ( t^{-\frac{7}{2}}). 
\end{aligned}
\end{equation}
Hence, we have
\begin{equation}
\begin{aligned}
E[X_{t}^2]= &\frac{g}{2 \sqrt{2} \pi ^{3/2} \sqrt{t}} + \frac{\sqrt{\frac{\pi }{2}}}{24 g t^{3/2}}-\frac{\pi ^{5/2}}{480 \sqrt{2} g^3 t^{5/2}} + \mathcal{O} ( t^{-\frac{7}{2}}),\quad t\rightarrow \infty.  \\
\end{aligned}
\end{equation}
The asymptotic expansion of variance as $t\rightarrow \infty$ is 
\begin{equation}
\begin{aligned}
  \mathrm{Var}[X_t]&=E[X_t^{2}]-(E[X_t])^{2}=\frac{g}{2 \sqrt{2} \pi ^{3/2} \sqrt{t}}- \frac{1}{4\pi t} + \frac{\sqrt{\frac{\pi }{2}}}{24 g t^{3/2}}-\frac{\pi ^{5/2}}{480 \sqrt{2} g^3 t^{5/2}} + \mathcal{O} ( t^{-\frac{7}{2}}).  \\
\end{aligned}
\end{equation}.

The  asymptotic expansion of higher moments of $E[A_{t}^{n}]$ can be derived by utilizing the recurrence relation provided in lemma \ref{thm:Matsumoto1}, 
\begin{equation}
\begin{aligned}
&E[A_{t}^{n}]= \frac{\Gamma (n)}{\sqrt{\pi } 2^{\frac{1}{2}-n}}t^{-\frac{1}{2}}+\mathcal{O} (t^{-1}), \quad  \frac{n}{2} \in \mathbb{Z}^{+}.\\
\end{aligned}
\end{equation}
where $\Gamma (x)$ is the gamma function. This relation also holds for a positive real number $n$ \cite{matsumoto2005exponential}. Then the leading order expansion of $E[X_{t}^n]$ is
\begin{equation}
\begin{aligned}
&E[X_{t}^{n}]=  (2 \pi )^{-\frac{n}{2}-\frac{1}{2}} g^{n-1} \Gamma \left(\frac{n}{2}\right) t^{-\frac{1}{2}} +\mathcal{O} (t^{-1}).
\end{aligned}
\end{equation}
It is easy to see that all normalized moment (e.g., skewness, kurtosis) diverges as $t\rightarrow \infty$, $\frac{E \left[ \left( X_{t}-E[X_{t}] \right)^{n} \right]}{\mathrm{var}(E_{t})^{n/2}} \sim (2 \pi )^{\frac{n-2}{4}} g^{\frac{n}{2}-1} \Gamma \left(\frac{n}{2}\right) t^{\frac{n-2}{4}}$, $n\geq 3$, meaning that the scalar becomes more and more intermittent during its decay. Whereas, the normalized moment in the random shear flow problem converges to a constant at long times \cite{bronski2000rigorous,bronski1997scalar}. Another example of divergent normalized moments is given in \cite{son1999turbulent}, where the flow is also white in time but has a more complicated spatial correlation. In that model, the kurtosis (flatness) grows as $t^{\frac{7}{4}}$.



\section{Conclusion and discussion}
\label{sec:conclusion}

In this paper, we considered an advection-diffusion equation (or a more general evolution equation) which involves a multiplicative Gaussian white noise.
We have derived the evolution equation for the ensemble average of the random scalar field by three different approaches, the Feynman-Kac formula, the limit of the Ornstein-Uhlenbeck process, and the calculus of white noise. All three different strategies yield the same evolution equation of ensemble-averaged scalar field, which demonstrates the validity of our results.   With our established evolution equation of ensemble average, we studied the passive scalar transport problem with two different types of flows, the periodic flow, and strain flow. For the periodic flow, by using the homogenization method, we show that the $N$-point correlation function of the random scalar flow is governed by a diffusion equation at long times, which generalized the conclusion in \cite{camassa2021persisting,ding2020ergodicity} from time-independent case to time-dependent case. For the strain flow, we have explicitly calculated the mean of the random scalar field and shown the statistic of this random field is related to the time integral of geometric Brownian motion.

  There are several future areas of exploration. First, in this paper, we focused on the system which is driven by the white noise process.  However, many applications concern the problem which involves a stochastic process with nonzero correlation. For example, the diffusing diffusivity model models the diffusivity as a function of OU process\cite{chechkin2017brownian,jain2017diffusing,tyagi2017non}. The study of dynamo \cite{bhat2015fluctuation,zel1984kinematic} and the passive-scalar decay problem\cite{vanneste2006intermittency,ding2021determinism,aiyer2017passive} concerns the renewing process. The modeling of Black-Scholes market \cite{elliott2003general}   involves fractional Brownian motion.  We expect that the strategies we presented in this paper could be generalized to problems that involve a more general stochastic process.
  Second, the standard homogenization method always yields an effective equation with constant coefficients. That could limit the valid range of parameters.
  Center manifold theory is a powerful tool to handle the scalar transport problem in bounded domain \cite{mercer1990centre}, and to derive the self-similarity solution of stochastic nonlinear diffusion reaction equation \cite{wang2013self}. A recent study \cite{ding2021determinism} has shown that, for the scalar transport problem in channel domain, the center manifold theory yields an effective equation with variable coefficients, which yields a long time limiting distribution of the random scalar field in the white noise case. Numerical simulation showed that this approximation performs better than the result from the standard homogenization method at an intermediate time scale. In addition, the effective equation suggests that a shear depending upon a Gaussian white noise process or a renewing process will induce a deterministic effective diffusivity at long time. We expect that the center manifold theory could yield a better approximation for the non-shear flow problem provided in equation \eqref{eq:advectionDiffusionNonDimension}. 
Last, we expect this averaged propagator to be useful in the analysis of fluctuations arising in many physical systems such as the ones alluded to in this paper.

\section{Acknowledgements}
We acknowledge funding received from the following National Science Foundation Grant Nos.:DMS-1910824; and Office of Naval Research  Grant No: ONR N00014-18-1-2490. Partial support for Lingyun Ding is gratefully acknowledged from the National Science Foundation, award NSF-DMS-1929298 from the Statistical and Applied Mathematical Sciences Institute. We
would also like to thank Andrew J. Majda for his comments and encouragement.

\subsection{Lists of abbreviations}
See table \ref{tab:abbreviations}.

\begin{table}[h]
 \centering
 \begin{tabular}{l|l}
\hline
  Full Form & Abbreviation\\
\hline\hline    
Ornstein-Uhlenbeck &OU\\
Partial differential equation & PDE \\
Probability distribution function & PDF\\
Stochastic differential equation &SDE\\   
\hline
 \end{tabular}
\caption{Lists of abbreviations.} \label{tab:abbreviations}
\end{table}
\bibliographystyle{elsarticle-harv}


\end{document}